\documentclass[10pt,twoside,a4paper]{article}
\usepackage{a4wide,times}
\usepackage{amsmath}
\usepackage{amsthm,amsfonts}
\usepackage{Latex-document}
\usepackage{amssymb}
\usepackage{graphicx}
\usepackage{float}
\usepackage{calc}
\usepackage{t1enc}
\usepackage[english]{babel}
\usepackage{1710}

\title{\bf Mathematical Modelling
  of \vskip -2mm the Cardiovascular System\vskip 6mm}
\author{A. Quarteroni\vspace*{-0.5cm}\thanks{MOX,
Department of Mathematics, Politecnico di Milano, Italy and
Institute of Mathematics (IMA), EPFL, Lausanne, Switzerland. E-mail: alfio.quarteroni@epfl.ch}}
\date{\vspace{-8mm}}

\begin{document}

\maketitle \thispagestyle{first} \setcounter{page}{839}

\begin{abstract}

\vskip 3mm

In this paper we will address the problem of developing mathematical models for
the numerical simulation of the human circulatory system. In
particular, we will focus our attention on the problem of
haemodynamics in large human arteries.
\vskip 4.5mm


\end{abstract}

\vskip 12mm

\section{Introduction} \label{section 1}\setzero
\vskip-5mm \hspace{5mm}

The simulation  of not only the physiological functioning of the blood
circulatory system, but also of
specific pathological circumstances is of utmost importance since cardiovascular diseases represent the leading cause
of death in developed countries, with a tremendous medical, social and
economic impact.

In the cardiovascular system, altered flow conditions, such as flow separation, flow reversal,
low and oscillatory shear stress areas, are recognised as important
factors in the development of arterial diseases.  A detailed
understanding of the local haemodynamics, the effect of vascular wall
modification on flow patterns and its long-term adaptation to surgical
procedures can have useful clinical applications. Some of these
phenomena are not well understood, making it difficult to foresee short
and long term evolution of the disease and the planning of the
therapeutic approach.  In this context, the mathematical models and
 numerical simulations can play a crucial role.

Blood flow interacts both mechanically and chemically with the vessel
walls. The mechanical
coupling requires algorithms that correctly describe the energy transfer
between the fluid (typically modelled by the Navier-Stokes equations)
and the structure of the vessel wall.   On the other hand, the flow
equations can be coupled with appropriate models that describe the
wall absorption of  bio-chemicals (e.g.\ oxygen, lipids, drugs,
etc.) and of their transport, diffusion and kinetics. Numerical simulations of this type may
help to understand the modifications in bio-chemical exchanges due to an
alteration of the flow field caused, for instance, by a stenosis (i.e.\
a localised narrowing of a vessel lumen, normally due to fat
accumulation).

The simulation of large and medium-size arteries is now sufficiently
advanced so to envisage the
applications of computer models to medical research and, in a medium
range,  to everyday medical practise.
For instance, simulating the flow in a coronary by-pass may
help  understanding the extent at which  its geometry influences the flow and in
turn the post-surgery evolution.  Also the study of the effects of a
vascular prosthesis  as well as the
study of artificial valve implants are areas which could benefit from  a sufficiently accurate simulation of blood flow field.

%
%
%
\par
In this paper we review the principal mathematical steps behind the
derivation of the coupled fluid-structure equations which model the
blood flow motion in large and medium-sized arteries. Then we mention
the way geometrical multiscale models, that combine mathematical models
set up in different spatial dimensions, can be conveniently used to
simulate the whole circulatory system.

\section{The coupled fluid-structure problem} \label{section 2}\setzero
\vskip-5mm \hspace{5mm}

\label{s:ale}

In this section we will treat the situation arising when the flow in a
vessel interacts mechanically with the wall structure. This aspect is
particularly relevant for blood flow in large arteries, where the vessel
wall radius may vary up to 10\percent because of the forces exerted by
the flowing blood stream.

We will first illustrate a framework for the Navier-Stokes equations in a
moving domain which is particularly convenient for the analysis and for
the set up of numerical solution methods.

\subsection{The Arbitrary Lagrangian Eulerian (ALE) formulation of the Navier-Stokes equations}

\vskip-5mm \hspace{5mm}

Navier-Stokes equations are usually derived according to the \emph{Eulerian} approach where
the independent spatial variables are the coordinates of a fixed
Eulerian system. When considering the flow inside a portion of a compliant
artery, we have  to compute the flow solution in a
\emph{computational domain} $\cdomain$ varying
with time.

\begin{figure}[hbp]
\begin{center}
\includegraphics[width=0.4\textwidth]{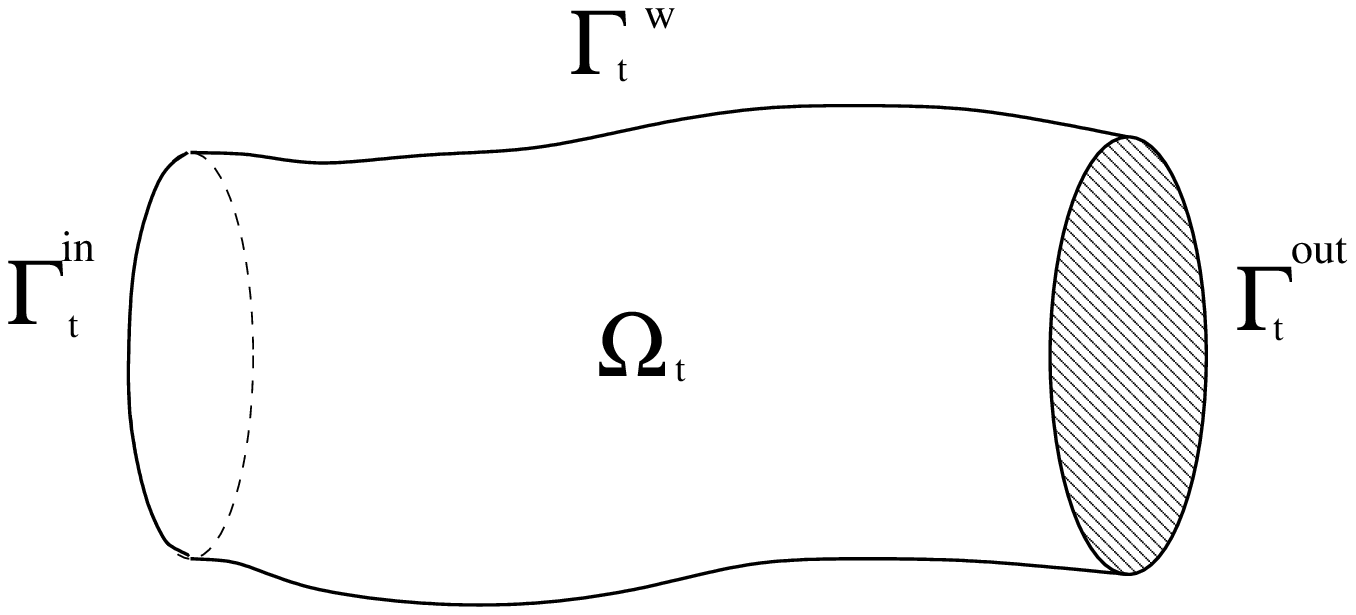}

\begin{minipage}[h]{10cm}
\caption{A simple model of a section of an artery. The vessel
  wall  $\wall_t$ is moving. The location
  along the $z$ axis of $\inflow$ and $\outflow$ are fixed.}
\label{f:computdomain}
\end{minipage}
\end{center}
\end{figure}

The boundary of $\cdomain$ may in general be subdivided into two parts.
The first part coincides with the physical fluid boundary, i.e.\ the
vessel wall $\wall_t$ (see Fig.~\ref{f:computdomain}),
which is moving under the effect of the flow
field. The other part of $\partial\cdomain$ corresponds to artificial
boundaries
 which delimit the
region of interest from the remaining part of the cardiovascular
system.

The ``artificial'' boundaries
are the inlet and outlet (or, using the medical terminology, the proximal
and distal) sections, here indicated by $\inflow$ and
$\outflow$, respectively.  The location of these boundaries is fixed a priori.
More precisely, $\inflow$ and $\outflow$ may change with time because of
the displacement of $\wall_t$, however they remain planar and their
position along the vessel axis is kept fixed. In this case the Eulerian approach becomes impractical.

A possible alternative would be to use the \emph{Lagrangian approach},
where we identify the computational domain on a reference configuration
$\crdomain$, e.g. that at the initial time $t=0$,
 and the corresponding domain in the current configuration
will be provided by the Lagrangian mapping
\begin{equation}
\cdomain=\ldomain=\LAG{t}(\crdomain),\,\,  t>0,
\end{equation}
which describes the motion of a material particle and whose time
derivative is the fluid velocity.
%
%
%
%
Since the fluid velocity at the wall is equal to the wall velocity, the
Lagrangian mapping effectively maps $\wall_0$ to the correct wall
position $\wall_t$ at each time $t$.  However, the artificial boundaries
in the reference configuration,  say $\Gamma^{i\!n}_0$ and
$\Gamma^{o\!u\!t}_0$,
will now be transported along the fluid trajectories. This is
unacceptable, particularly for
a relatively large time interval as $\cdomain$ rapidly becomes
highly distorted.

A more convenient situation is the one when, even if the wall is moving, one keeps
the inlet and outlet boundaries at the same spatial location along the
vessel axis. With that purpose, we introduce the \emph{Arbitrary Lagrangian
Eulerian} (ALE) mapping
$\ALE: \crdomain\rightarrow\adomain, \qquad \Y \rightarrow \y(t,\Y)=\ALE(\Y)$,
which provides the spatial coordinates $(t,\y)$ in terms of the
so-called \emph{ALE coordinates} $(t,\Y)$, with the basic requirement that
$\ALE$ retrieves, at each time $t>0$, the desired
computational domain, $\cdomain=\adomain=\ALE(\crdomain), \, t\ge 0$.

The ALE mapping should be continuous and bijective in
$\overline{\crdomain}$. Once given, we may define the \emph{domain}
velocity field as
\begin{equation}
{\mvel}=\frac{\partial\ALE}{\partial t}\circ\ALEI,
\end{equation}
where the composition operator applies only to the spatial coordinates.
The ALE time derivative of a function $f: \timeint\times\ctdomain\rightarrow\R$, which we denote by
$\aleder{f}$, is defined as

\begin{equation}
\label{e:matderoo}
\aleder{f}: \timeint\times\ctdomain\rightarrow \R, \qquad \aleder f=
\timeder{\widetilde{f}}\circ\ALEI.
\end{equation}
This definition is readily extended to vector valued functions. The ALE
derivative is related to the Eulerian (partial) time derivative by the relation

\begin{equation}
\label{e:aleder}
\boxed{\aleder{f}=\eulerder{f} + \mvel \Dot \Grad{f},}
\end{equation}
where the gradient is made with respect to the $\y$-coordinates.
\smallskip

The Navier-Stokes equations may be formulated in
order to put into evidence the ALE time derivative, obtaining
\begin{equation}
\label{e:anavsto}
\begin{array}{l}
\aleder{\vel} +  \aaconvop\vel + \Grad{p}
-\DIV{\mathbf{T}(\vel)}=\fbody, \\[5.mm]
\Div{\vel}=0,
\end{array}\quad \text{in } \ctdomain, \, t>0,
\end{equation}
where $\mathbf{T}$ is the Cauchy stress tensor, which for Newtonian fluid is
given by  $\mathbf{T}=\viscok(\nabla\mathbf{u}+\nabla\mathbf{u}^T)$,
being $\viscok$ the blood kinematic viscosity.

When considering small vessels, accounting for non-Newtonian behaviour
of blood becomes crucial. In that case the functional dependence of
$\mathbf{T}$ on $\mathbf{u}$ becomes more complex, see for instance
\cite{rajagopal.kr:mechanics}.

\subsection{The structure model}
\label{s:structure}

\vskip-5mm \hspace{5mm}

The vascular wall has a very complex nature and devising an accurate
model for its mechanical behaviour is rather difficult. Its structure is
indeed formed by many layers with different mechanical
characteristics \cite{fung:biomechanics,holzapfel00}, which are usually
in a pre-stressed state.
Moreover, experimental results obtained by
specimens are only partially significant. Indeed, the vascular wall is a
living tissue with the presence of muscular cells which contribute to
its mechanical behaviour.  It may then be expected that the dead tissue
used in the laboratory will have different mechanical characteristics
than the living one. Moreover, the arterial mechanics depend also on the
type of the surrounding tissues, an aspect almost impossible to
reproduce in a laboratory. As we have already pointed out, the displacements cannot be considered small (at least in large
arteries where  the radius may vary up to a few percent during the
systolic phase). Consequently, an
appropriate model for the structure displacement $\boldsymbol{\eta}$ reads
\[
\rho_w\frac{\partial^2\boldsymbol{\eta}}{\partial
  t^2}-\boldsymbol{\operatorname{div}}\boldsymbol{\sigma}(\boldsymbol{\eta},\frac{\partial \boldsymbol{\eta}}{\partial t})
=\mathbf{f},\quad \text{in }\Omega_t^s,\, t>0,
\]
where $\Omega_t^s$
indicates the current configuration and $\boldsymbol{\sigma}$ is the
Cauchy stress tensor. The latter may depend
on the structure velocity because of viscoelasticity. A full Lagrange
formulation for the structure on a fixed reference configuration
$\Omega_0^s$ may be obtained by the usual Lagrange and Piola transformation (see, e.g., \cite{ciarlet:mathematical}).
A  general framework to derive constitutive equations for arterial walls
is reported in \cite{holzapfel00}.

It is the
role of mathematical modelling to find reasonable simplifying
assumptions by which major physical characteristics remain present,
yet the problem becomes computationally actractive.
In particular, a simpler model may be obtained by considering only
displacements in the radial direction and a cylindrical geometry for the
vessel. Furthermore if we neglect the geometrical non-linearities (which
correspond to assume small
displacements), as well as the variations along the radial
directions (small thickness assumption)  we obtain the following
``generalised string model''
\cite{quarteroni.a.tuveri.m.ea:computational,QF02:modelling} for the
evolution of the radial displacement $\eta=R-R_0$,
\begin{equation}
\label{e:generstring}
\boxed{
\frac{\partial^2 \vdisp}{\partial t^2} -a \frac{\partial^2
 \vdisp}{\partial z^2}+ b\vdisp-c\frac{\partial^3 \vdisp}{\partial
t\partial z^2}=\f,\quad \text{in
} \vesselref,\, t>0, }
\end{equation}
where $a>0$ and  $b>0$ are parameters linked to the vessel geometry and
mechanical characteristics, $c>0$ is a viscoelastic parameter and  $H$ is
a forcing term which depends on the action of the fluid, as we will see
in (\ref{eq:interaction}).
More
details as well as the derivation of this model are found in the cited references.

Here $\vesselref$ is the reference configuration for the structure
\begin{equation*}
\vesselref=\{(r,\theta,z):\,r=\radiusref(z),\, \theta\in[0,2\pi),\, z\in[0,L]\},
\end{equation*}
where $L$ indicates the length of the arterial element under
consideration.  In our cylindrical coordinate system $(r,\theta,z)$, the
$z$ coordinate is aligned along the vessel axes and a plane
$z=\overline{z}\,(=\text{constant})$ defines an \emph{axial section}.
%

\subsection{Coupling with the structure model}
\label{s:alecouple}

\vskip-5mm \hspace{5mm}

We now study the properties of the coupled fluid-structure problem,
using for the structure the generalised string model
(\ref{e:generstring}).  We will  take $\normal$ always to be
the outwardly vector normal to the fluid domain boundary. Furthermore we
define $g$ as the metric function so that the elemental surface measure
$d\sigma$ on $\vessel$ is related to the corresponding measure
$d\sigma_0$ on $\vesselref$ by  $d\sigma= g\, d\sigma_0$.

We will then address the following problem:
\emph{
For all $t>0$, find $\vel$, $\press$ , $\vdisp$ such that
}
\begin{equation}
\label{eq:interaction}
\begin{cases}
\vel,\pres \quad \text{satisfy problem (\ref{e:anavsto})},\\
\vdisp \quad \text{satisfies problem (\ref{e:generstring})},\\
\vel\circ\ALE=\Frac{\partial\vdisp}{\partial t}\mathbf{e}_r,\qquad
\text{on } \vesselref,\\
\f= g\displaystyle\frac{\rho}{\rho_w h_0}
\left[(\press-\press_0)\mathbf{n}-\mathbf{T}(\mathbf{u})\cdot\mathbf{n}\right]\cdot\mathbf{e}_r\qquad
\text{on } \vesselref.
\end{cases}
\end{equation}
Here, $p_0$ is the pressure acting at the exterior of the vessel,
$\mathbf{e}_r$ is the radial unit vector, $\rho_w$ and $\rho$ are the
wall and fluid densities, respectively, while $h_0$ is the wall
thickness. The system is complemented
by appropriate boundary and initial conditions.

We may then recognise the sources of the coupling between the fluid and
the structure models, which are twofold. In view of a possible
iterative solution strategy, the fluid solution
provides the value of $H$, which is function of the fluid stresses at
the wall. On the other hand, the movement of the vessel wall
  modifies the geometry on which the fluid equations must be solved,
  besides providing Dirichlet boundary conditions for the fluid velocity in
  correspondence to the vessel wall.

\noindent{\bf Remark 2.1.}
  We may note that the non-linear convective term in the Navier-Stokes
  equations is crucial to obtain the well-posedness of the coupled
  problem, because it
  generates a boundary term which compensates that coming from the
  treatment of the acceleration term. These two contributions are
  indeed only present in the case of a moving boundary. See
  \cite{QF02:modelling} and \cite{veiga02}.

\subsection{Numerical solution of the coupled fluid-structure problem}\label{s:coupling}

\vskip-5mm \hspace{5mm}

In this section we describe an algorithm that at each
time-level allows the decoupling of the sub-problem related to the
fluid from that  related to the vessel wall.  As usual, $t^k$,
$k=0,1,\ldots$  denotes the k-th discrete time level; $\Delta t >0$
is the time-step, while $v^k$ is the approximation of the function
(either  scalar or vector) $v$ at time $t^k$.

The numerical solution of the fluid-structure interaction problem
(\ref{eq:interaction}) will be carried out by
constructing a proper finite element approximation of each
sub-problem.  In particular, for the fluid we need to devise a finite
element formulation suitable for moving domains (or, more precisely,
moving grids). In this respect, the ALE formulation will
provide an appropriate framework.

\begin{figure}[hbp]
\begin{center}
\includegraphics[width=0.45\textwidth,height=2.5cm]{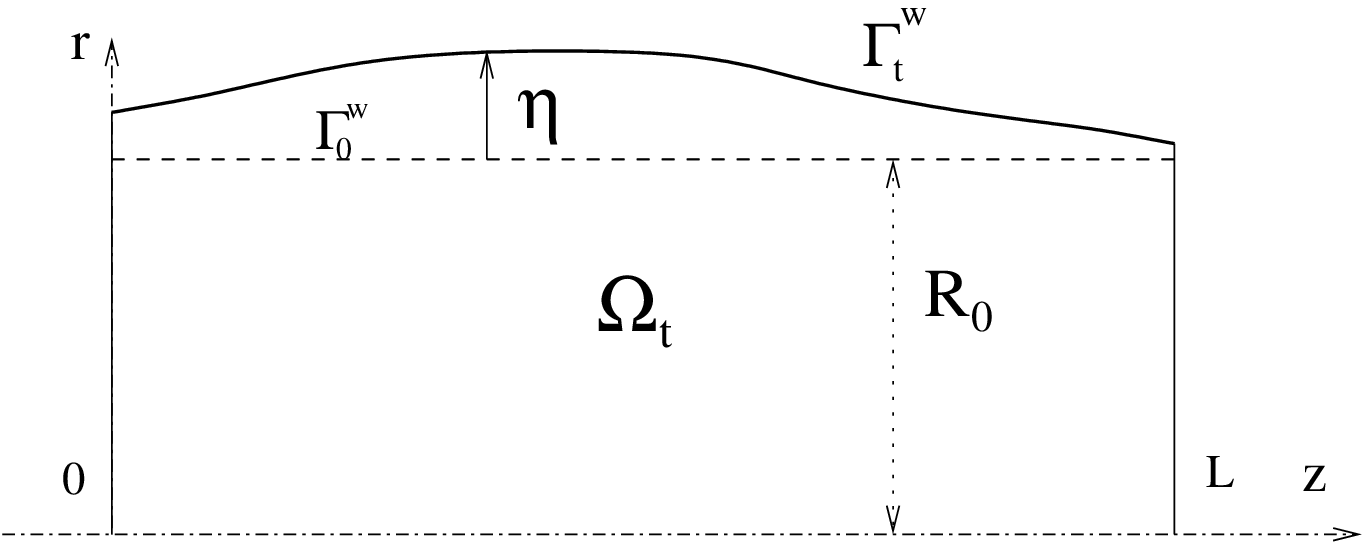}
\hfill\includegraphics[width=0.45\textwidth,height=2.5cm]{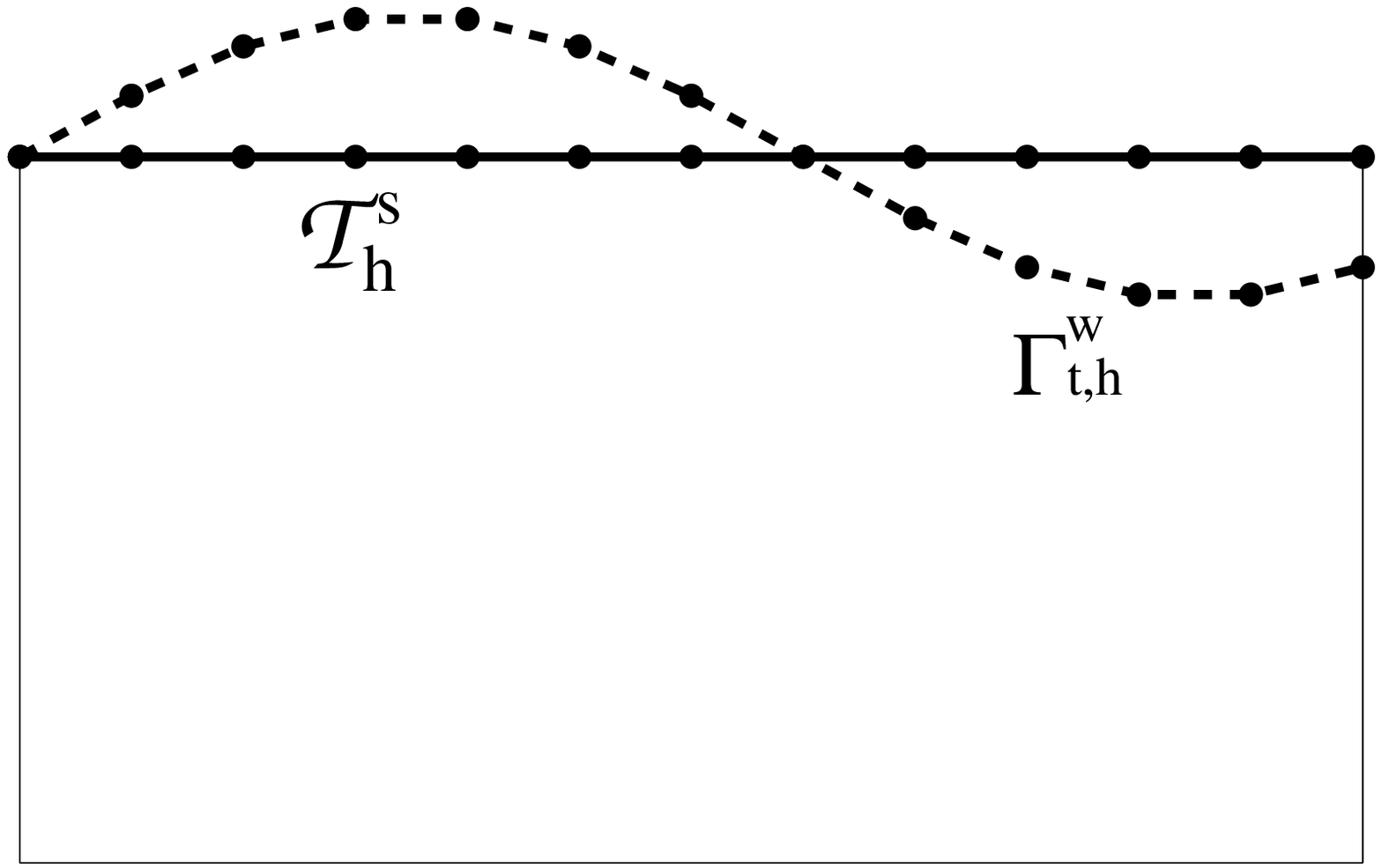}

\begin{minipage}[h]{10cm}
\caption{A simple fluid-structure interaction problem. On the left
the  domain definition and on the right the discretized vessel
wall corresponding to a possible value of
  $\vdisp_h$.}
\label{fig:twodaxi}
\end{minipage}
\end{center}
\end{figure}

To better illustrate the situation we refer to Fig.~\ref{fig:twodaxi}
where we have drawn a 2D fluid structure interaction problem (only the
upper portion of the vessel is reported). For the sake of simplicity we have considered only  a
two-dimensional fluid structure problem, yet
the algorithm here presented may be readily extended to more complex
situations and three dimensional problems.

\begin{figure}[hbt]
\begin{center}
\includegraphics[width=0.85\textwidth,height=3.8cm]{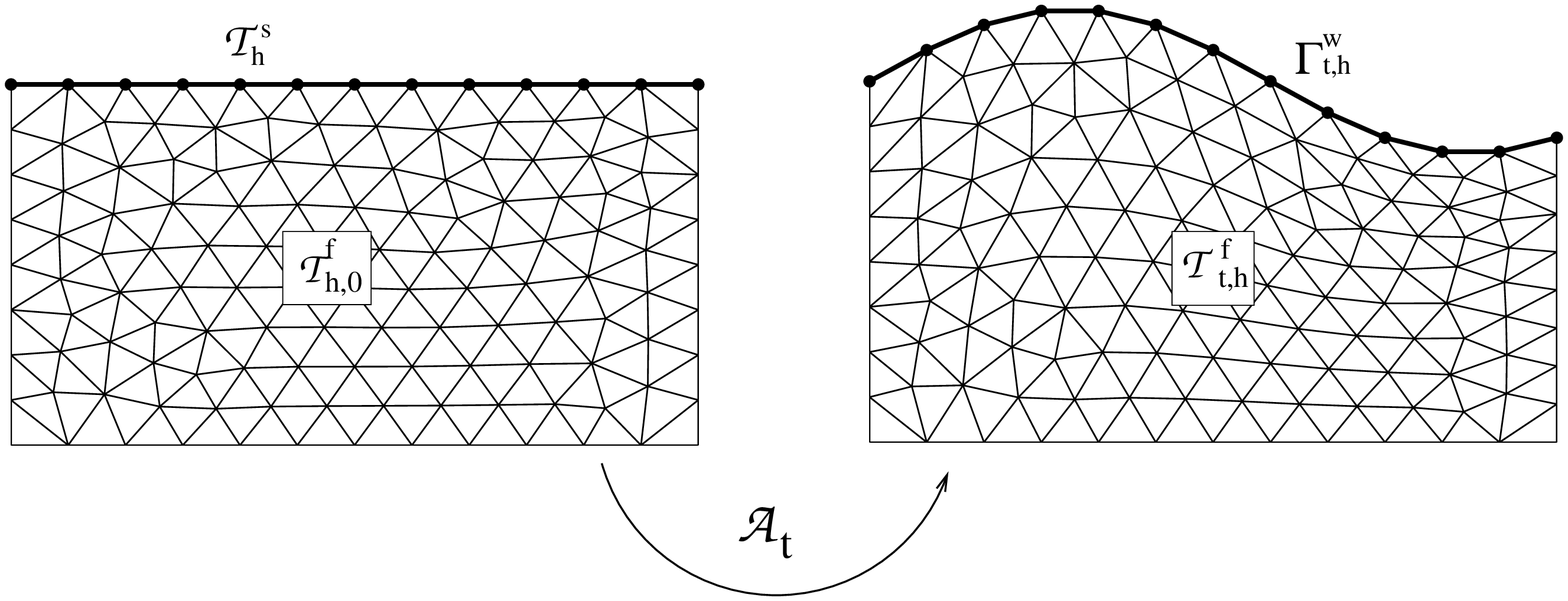}

\begin{minipage}[h]{10cm}
\caption{The triangulation used for the fluid problem at each time $t$
  is the the image through a map $\ALE$ of a mesh constructed on $\Omega_0$.}
\label{fig:disc_ale}
\end{minipage}
\end{center}
\end{figure}

The structure on $\vesselref$ will be discretized by means of a grid
 $\vesselmesh$ and  employing piece-wise linear continuous (P1)
finite elements  to represent the approximate vessel wall displacement
$\vdisp_h$.  The position at
time $t$ of the discretized vessel wall boundary,
corresponding to the discrete displacement field $\vdisp_h(t)$, is
indicated by $\vesseldiscrete$.
Consequently, the fluid
domain will be represented at every time by a polygon, which we indicate
by $\Omega_{t,h}$. Its triangulation $\fluidmesh$ will be constructed as
the image by an appropriate ALE mapping $\ALE$ of a triangulation
$\fluidmeshref$ of $\Omega_0$, as shown in Fig.~\ref{fig:disc_ale}.
Correspondingly, $\Omega_{t,h}=\ALE{\Omega_{0,h}}$, where $\Omega_{0,h}$
is the discretisation of $\Omega_0$ induced by $\fluidmeshref$.
The trace of $\fluidmeshref$ on $\vesselref$ will coincide with the
$\vesselmesh$ of the vessel wall, i.e. we here consider
\emph{geometrically conforming} finite elements between the fluid and
the structure.  The possibility of using  geometrically non-conforming
finite elements has been investigated in
\cite{grandmont98:_noncon}.  \smallskip

Having at disposal the discrete displacement field $\vdisp^{k+1}_h$ at
$t=t^{k+1}$ and thus the position of the domain boundary $\partial
\Omega_{t^{k+1},h}$,  the set up of  a map ${\cal A}_{t^{k+1}}$ such that
${\cal A}_{t^{k+1}}(\fluidmeshref)$ is an acceptable finite element mesh
for the fluid domain is not a simple task. However, if we
can assume that $\Omega_{t,h}$ is convex for all $t$ and that the
displacements are relatively small, the following technique, known as
\emph{harmonic extension}, may well serve the purpose. If  $\mathbf{X}_h$
indicates  the P1 finite element vector space associated to $\fluidmeshref$
and  $\mathbf{g}_h:\partial \Omega_{0,h}\rightarrow \partial\Omega_{t^{k+1},h}$
is the
function describing the fluid domain boundary, we build the map by
seeking $\mathbf{y}_h\in \mathbf{X}_h$ such that
\begin{equation}
\label{eq:map}
\displaystyle\int_{\Omega_0} \Grad{\mathbf{y}_h}:\Grad{\mathbf{z}_h}=0\quad \forall
\mathbf{z}_h\in \mathbf{X}_h^0,\qquad
\mathbf{y}_h=\mathbf{g}_h,\quad \text{on }\, \partial\Omega_{0,h},
\end{equation}
and then setting ${\cal A}_{t^{k+1}}(\Y)=\mathbf{y}_h(\Y),\quad
\forall\Y\in\Omega_{0,h}$. A more general discussion on the
construction of the ALE mapping may be found in
\cite{formaggia.nobile:stability,nobile01:_numer_approx} as well
as in \cite{gastaldi01:_lagran_euler}.

\noindent {\bf Remark 2.2.}
  Adopting P1 elements for the construction of the ALE map ensures that the triangles
  of ${\cal T}^f_{h,0}$ are mapped into triangles, thus ${\cal T}^f_{h,t}$
  is a valid triangulation, under the requirement  of invertibility
  of the map (which is assured if the domain is convex and the
  wall displacements are small).

As for the time evolution, we may adopt a linear time variation within
each time slab $[t^{k},t^{k+1}]$ by setting
\[
\ALE=\frac{t-t^{k}}{\Delta t}{\cal A}_{t^{k+1}}-\frac{t-t^{k+1}}{\Delta
  t}{\cal A}_{t^{k}}, \quad t\in[t^{k},t^{k+1}].
\]
Then, the corresponding domain velocity $\mathbf{w}_h$ will be constant
on each time slab.

\subsubsection{The iterative algorithm}

\vskip-5mm \hspace{5mm}

We are now in the position of describing an iteration algorithm for
the solution of the coupled problem.  As usual, we assume that all
quantities are available at $t=t^k$, $k\ge0$, provided either by previous calculations
or by the initial data and we wish to advance to the new time step
$t^{k+1}$.  For ease of notation we here omit the subscript $h$, with the
understanding that we are referring exclusively to finite element
quantities.

The algorithm requires to choose a \emph{tolerance} $\tau>0$, which is
used to test the convergence of the procedure, and a \emph{relaxation
  parameter} $0<\theta\le 1$.  In what follows, the subscript $j\ge 0$ denotes the sub-iteration counter.

The algorithm reads:
\renewcommand{\labelenumi}{A\arabic{enumi}}
\renewcommand{\labelenumii}{A\arabic{enumi}.\arabic{enumii}}
\begin{enumerate}
\setcounter{enumi}{0}
\item Extrapolate the vessel wall structure displacements  and velocity:
\[
\vdisp_{(0)}^{k+1}=\vdisp^{k}+\Delta t \vdispv^{k},\qquad \vdispv_{(0)}^{k+1}=\vdispv^{k}.
\]
\item Set $j=0$.
\begin{enumerate}
\item By using $\vdisp_{(j)}^{k+1}$ compute
  the new grid for the fluid domain $\Omega_{t}$  and the ALE map by
  solving the harmonic extension problem (\ref{eq:map}).\label{item:a}
\item Approximate the Navier-Stokes problem to compute
  $\vel^{k+1}_{(j+1)}$ and $\press^{k+1}_{(j+1)}$, using as
  velocity on the wall boundary the one calculated from
  $\vdispv_{(j)}^{k+1}$.
\item Approximate the structure problem to compute $\vdisp_{*}^{k+1}$ and
  $\vdispv_{*}^{k+1}$ using
 $\vel^{k+1}_{(j+1)}$ and $\press^{k+1}_{(j+1)}$ to recover the
forcing term $\f$.
\item Unless
$\Vert
\vdisp_{*}^{k+1}-\vdisp_{(j)}^{k+1}\Vert_{L^2(\vesselref)}+
\Vert\vdispv_{*}^{k+1}-\vdispv_{(j)}^{k+1}\Vert_{L^2(\vesselref)}
\le\tau$, set
\[
\vdisp_{(j+1)}^{k+1}=\theta \vdisp_{(j)}^{k+1}+
(1-\theta)\vdisp_{*}^{k+1},\,\, \vdispv_{(j+1)}^{k+1}=\theta \vdispv_{(j)}^{k+1}+ (1-\theta)\vdispv_{*}^{k+1},
\]
and $j\leftarrow j+1$. Then return to
step~\ref{item:a}.
\end{enumerate}
\item Set
\[
\vdisp^{k+1}=\vdisp_{*}^{k+1}, \qquad \vdispv^{k+1}=\vdispv_{*}^{k+1}.
\]
\[
\vel^{k+1}=\vel_{(j+1)}^{k+1}, \qquad \press^{k+1}=\press_{(j+1)}^{k+1}.
\]
\end{enumerate}

If the algorithm converges, then
$\lim_{j\rightarrow\infty}\vel^{k+1}_{(j)}=\vel^{k+1}$ and
$\lim_{j\rightarrow\infty}\vdisp^{k+1}_{(j)}=\vdisp^{k+1}$, where
$\vel^{k+1}$ and $\vdisp^{k+1}$ are the approximate solution of the
coupled problem at time step $t^{k+1}$.

The algorithm entails, at each sub-iteration, the computation of the
 equation for the structure mechanics, the Navier-Stokes equations and the
solution of two Laplace equations (\ref{eq:map}), one for every
displacement component. It is therefore quite computationally
expensive. Alternatively, less implicit formulations may be adopted, see
for instance \cite{piperno01:_partit}, yet it has been found that for
the problem at hand a strong coupling between fluid and structure must
be maintained also at discrete level in order to have stable algorithms \cite{nobile01:_numer_approx}.

\section{Multiscale modelling of the cardiovascular system}\label{section 4}\setzero
\vskip-5mm \hspace{5mm}

\begin{figure}[p]
\begin{center}

\includegraphics[width=0.75\textwidth]{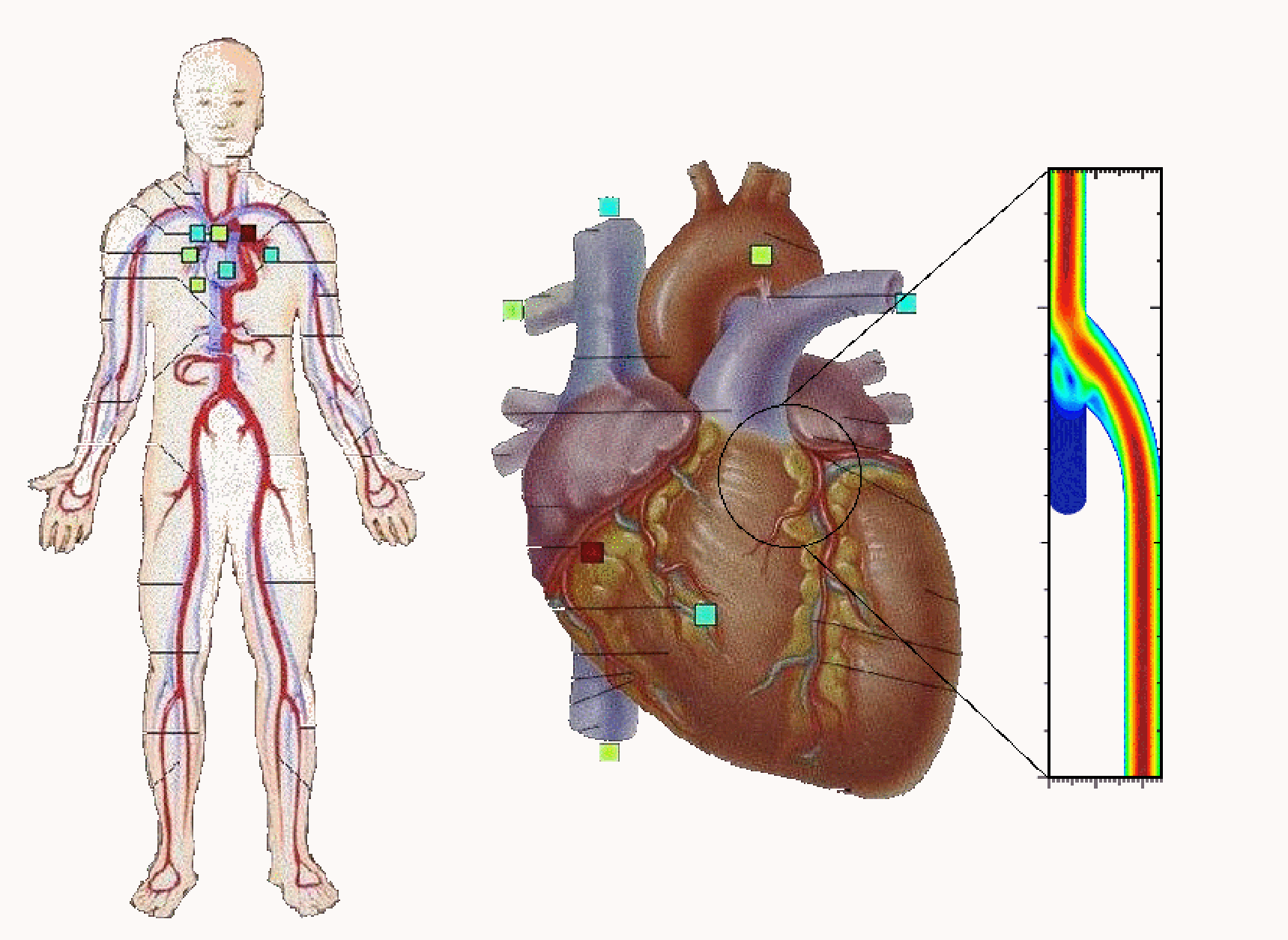}\\
\includegraphics[width=0.75\textwidth]{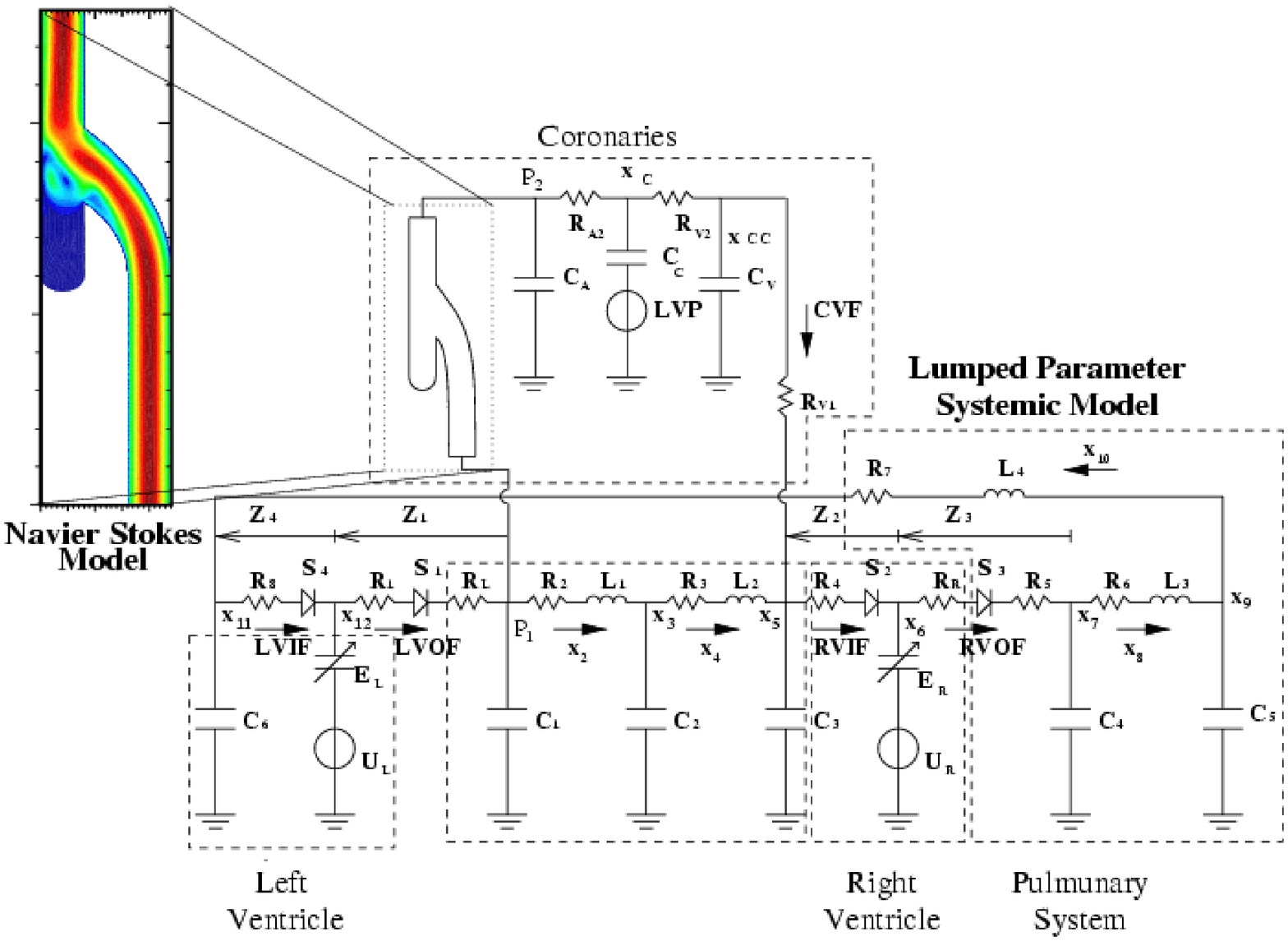}

\begin{minipage}[h]{10cm}
\caption{ An example of  multiscale simulation of
  blood flow, with the interplay between three-dimensional,
  one-dimensional and lumped parameters models. On top we show a global
  model of the circulatory system where a coronary by-pass is being
  simulated by a Navier-Stokes fluid-structure interaction model.
  The rest of the circulatory system is described by means of a lumped
  parameter model, based
  on the solution of a system of ODEs, is here represented by an
  electrical circuit analog in the bottom part of the figure.
 }\label{fig:allcoupled}
\end{minipage}
\end{center}
\end{figure}

The
cardiovascular system is highly integrated. In many cases, to isolate
the part of interest from the rest of the system would require
specification of point-wise boundary data on artificial boundary
sections. These are difficult to pre-determine.  To account for the effect of the global
circulatory system when focusing on specific regions we propose to
integrate a hierarchy of models operating at different ``scales''.  At
the highest level we have the full three-dimensional fluid-structure
interaction problem. This will be used where details of local flow
fields are needed. At the lowest level, we use lumped parameter models
based on the resolution of systems of non-linear ordinary
algebraic-differential equations  for averaged mass flow and
pressure. The latter models are often described by help of an analogy with
an electrical circuit, where the voltage represents blood pressure and
the current the flow rate. They are  well fit to supply the more sophisticated models
with  the effects of the circulation in small vessels,
the capillary bed, the venous system, as well as the action of the heart.
A transition between the two extrema could be achieved by convenient
one-dimensional models expressed by a first order non-linear
hyperbolic system (see Fig.  \ref{fig:allcoupled}). The derivation of the one-dimensional model and
a possible  numerical implementation may be found in \cite{formaggia02,formaggia02:_one,QF02:modelling}.

An analysis of the coupling between fluid-structure models and one
dimensional models may be found in
\cite{formaggia.gerbeau.eal:coupling}, while the direct coupling by
lumped parameter models and Navier-Stokes model is found in
\cite{quarteroni01:_coupl,quarteroni02:_analy_geomet_multis_model_based},
the coupling between lumped parameter models and one-dimensional models
is also treated in \cite{formaggia.nobile.ea:multiscale}.

\noindent{\bf Acknowledgements.} The research activity here
described has been supported by various Swiss and Italian research
agencies and isntitutions, in particular the Swiss National
Research Fund (FNS), and the Italian CNR, Ministry of Education
(MIUR). The author also thanks Luca Formaggia for his contribution
to the preparation of this paper.

\label{lastpage}

\end{document}